\documentclass[pdflatex,sn-mathphys-num]{sn-jnl}% Math and Physical Sciences Numbered Reference Style

\usepackage{graphicx}%
\usepackage{multirow}%
\usepackage{amsmath,amssymb,amsfonts}%
\usepackage{amsthm}%
\usepackage{mathrsfs}%
\usepackage[title]{appendix}%
\usepackage{xcolor}%
\usepackage{textcomp}%
\usepackage{manyfoot}%
\usepackage{booktabs}%
\usepackage{algorithm}%
\usepackage{algorithmicx}%
\usepackage{algpseudocode}%
\usepackage{listings}%

\hypersetup{
	colorlinks = true,
	urlcolor   = blue,
	citecolor  = black,
}
\numberwithin{equation}{section}
\theoremstyle{thmstyleone}
\newtheorem{theorem}{Theorem}
\newtheorem{proposition}[theorem]{Proposition}
\newtheorem{lemma}[theorem]{Lemma}
\newtheorem{corollary}[theorem]{Corollary}

\newcommand{\R}{\mathbb{R}}
\newcommand{\abs}[2][]{#1\lvert #2 #1\rvert}

\newcommand{\ucrest}{\hat u} % speed at the crest
\newcommand{\uni}{\mathrm{uni}} % subscript for unidirectional waves

\begin{document}

\title{On the vorticity threshold for steady water waves}

\author*[1]{\fnm{Evgeniy} \sur{Lokharu}}\email{evgeniy.lokharu@math.lu.se}

\author[2]{\fnm{Miles H.} \sur{Wheeler}}\email{mw2319@bath.ac.uk}

\affil*[1]{\orgdiv{Centre for Mathematical Sciences}, \orgname{Lund University}, \orgaddress{\city{Lund}, \postcode{SE-223 62}, \country{Sweden}}}

\affil[2]{\orgdiv{Department of Mathematical Sciences}, \orgname{University of Bath}, \orgaddress{\city{Bath}, \postcode{BA2 7AY}, \country{UK}}}

\abstract{
  One of the most striking features of two-dimensional steady water waves with adverse vorticity is the emergence of stagnation points in the flow, even along branches of solutions that are initially unidirectional. In this paper we show that, for constant adverse vorticity above a certain explicit threshold, the first stagnation point must appear on the bed directly below the crest. This is in contrast to waves with favourable or even small adverse vorticity, where it is known that bottom stagnation cannot occur, leading to extreme waves exhibiting surface singularities. 
  % Miles: what's the reference for bottom stagnation not happening with small adverse vorticity?
  We also establish several related inequalities for unidirectional waves with strong adverse vorticity, including a new upper bound on the amplitude.
}

\maketitle

\section{Introduction} \label{s:introduction}

The present paper concerns steady water waves with large vorticity. We assume an ideal fluid of constant unit density occupies a two-dimensional region of finite depth, bounded below by a flat bottom and above by a free surface. Neglecting the effects of surface tension and air motion, the hydrodynamic pressure is continuous across this interface and coincides with the constant air pressure there. We work with the full steady Euler equations, formulated in non-dimensional variables where the relative mass flux and the gravity are scaled to unity; see Section~\ref{s:formulation} for a detailed description. The (non-dimensional) relative velocity field is denoted by $(u,v)$, and the surface is always given by the graph $y = \eta(x) > 0$. The flow is assumed to be rotational, where the vorticity function $\omega = v_x-u_y$ is constant and \emph{adverse}, that is $\omega = -\gamma < 0$. While even small waves with constant adverse vorticity may have counter currents \cite{Wahlen09,Kozlov2020}, we are focused here on waves that are \emph{unidirectional} in the sense that $u > 0$ everywhere in the fluid.

We consider only smooth classical solutions, either Stokes waves or solitary waves. Here a Stokes wave is one with a periodic profile having exactly one crest and trough in each minimal period, while a solitary wave is one that converges uniformly to a flat background state at infinity. In the unidirectional setting, it is known that both Stokes \cite{ConstantinEscher04b,ConstantinEhrnstromWahlen06,Hur07} and solitary \cite{Hur08, Kozlov2021} waves are necessarily symmetric, and that solitary waves are `waves of elevation' with exactly one crest. Moreover, solitary waves have a supercritical wave speed \cite{Keady1978, Wheeler15, Kozlov2021}. By a laminar flow we mean a trivial solution where the velocity is purely horizontal and the free surface is completely flat.

Our first main result can now be stated as follows.

\begin{theorem}\label{thm:main}
  Along any continuous family of unidirectional Stokes or solitary waves, connected to a laminar flow and with constant adverse (dimensionless) vorticity
  \begin{align*}
    -\omega = \gamma > \gamma_\star\approx 1.37184,
  \end{align*}
  the horizontal velocity satisfies the lower bound
  \begin{equation} \label{thm:main:psiy}
    u \geq \tfrac12 \gamma \eta 
    \quad
    \textup{on the surface $y=\eta(x)$}.
  \end{equation}
  In particular, extreme or overhanging waves cannot form before a stagnation point has emerged on the bottom.
\end{theorem}

This result is consistent with numerical observations of a certain threshold value of adverse vorticity. Below the threshold value waves approach surface stagnation as in the irrotational setting, while for adverse constant vorticity above the threshold one first observes bottom stagnation and the formation of eddies after that. Such a transition was reported by \cite{Teles-daSilvaPeregrine88} for solitary waves, and by \cite{Ko08large,Ko08effect} for periodic waves. More recently, \cite{Ribeiro17} numerically calculated the regions in parameter space where the number of stagnation points changes for periodic waves. We refer the reader to the introduction of the latter reference for a more thorough discussion.

% Miles: My best guess is that Ko & Strauss see a transition at \gamma \approx 0.45. The solitary paper is harder because they do not report the flux, but doing some rough calculations looking at their figure 11 I'm getting that their smallest value of \gamma with eddies has to be at least 0.64. Not sure any of this is worth mentioning, though.

By combining Theorem~\ref{thm:main} with existing global bifurcation results \cite{Haziot2023}, we are able to obtain the existence of solitary waves with a single stagnation point. Similar statements can be shown for Stokes waves.

\begin{theorem}\label{thm:main2}
  For any $\gamma > \gamma_\star$, there exists a smooth solitary wave for which $u > 0$ everywhere in the fluid except one point at the bottom beneath the crest, where $u = 0$.
\end{theorem}

Asymptotically sharper estimates than \eqref{thm:main:psiy} can be obtained for any large adverse vorticity. To simplify the formulas we introduce the shorthand
\begin{align*}
  \hat \eta = \sup_{x \in \R} \eta(x),
  \qquad 
  \check \eta = \inf_{x \in \R} \eta(x)
\end{align*}
for the depth of the fluid at a crest or trough, with the convention that the `trough' of a solitary wave of elevation is located at infinity.

\begin{theorem}\label{thm:main1}
  Along any curve of unidirectional solitary waves, with constant adverse vorticity $\gamma>\gamma_\star$ and bifurcating from a shear flow, the following is true.
  \begin{itemize}
  \item[(i)] Everywhere along the surface we have
    \[
      u \geq \gamma \check{\eta} - O(\gamma^{-5/2}) = \sqrt{2\gamma} - O(\gamma^{-1}).
    \]
  \item[(ii)] At every crest we have
    \[
      u \geq \sqrt{2\gamma} - O\big(\gamma^{-\tfrac52}\big)
      \quad \textup{and} \quad 
      \hat{\eta} \leq \sqrt{\frac{2}{\gamma}} + O\big(\gamma^{-\tfrac72}\big).
    \]
  \item[(iii)] The amplitude satisfies
    \[	
    \hat{\eta} - \check{\eta} = O(\gamma^{-2}).
    \]
  \end{itemize}
  The same is true for unidirectional Stokes waves provided the Bernoulli constant $r$ satisfies the bound
  \begin{equation} \label{eqn:rbound}
    r \le \gamma + \sqrt{\frac 2\gamma}.
  \end{equation}
  Here the $O(\,\cdot\,)$ terms refer to the limit as $\gamma \to \infty$, and the implied constants are absolute. In particular, they are independent of the solution and of the vorticity.
\end{theorem}
The bound \eqref{eqn:rbound} automatically holds, for instance, for curves of solutions with fixed $r,\gamma$ and variable wavelength. See Section~\ref{s:preliminaries} for more details.

At first glance it may seem surprising that strong adverse vorticity has the same sort of damping effect on amplitude as favorable vorticity \cite{Constantin2021, Lokharu2023, Haziot2024}, especially when numerical observations predict that a strong adverse current should \emph{increase} wave height and steepness \cite{Teles-daSilvaPeregrine88, SwanCummingJames01,Vanden-Broeck94}. The explanation for this discrepancy is that such waves are only possible after the unidirectionality assumption has been violated, highlighting the importance of critical layers for the existence of large-amplitude solutions.

\section{Problem formulation} \label{s:formulation}

Working in a reference frame traveling with the wave, the flow inside the fluid is 
modeled by the steady Euler equations
\begin{subequations}\label{eqn:trav}
  \begin{align}
    \label{eqn:u}
    u u_x + vu_y & = -P_x,   \\
    \label{eqn:v}
    u v_x + vv_y & = -P_y-g, \\
    \label{eqn:incomp}
    u_x + v_y &= 0.
  \end{align}
  Here $(u,v)$ is the two-dimensional relative velocity field, $P$ is the pressure, $g > 0$ is the gravitational constant, and the equations hold true within the two-dimensional fluid domain $D = \{ (x,y): x\in \R,\, 0 < y < \eta(x)\}$. The kinematic and dynamic boundary conditions are
  \begin{alignat}{2}
    \label{eqn:kinbot}
    v &= 0&\qquad& \text{on } y=0,\\
    \label{eqn:kintop}
    v &= u\eta_x && \text{on } y=\eta,\\
    \label{eqn:dyn}
    P &= P_{\mathrm{atm}} && \text{on } y=\eta,
  \end{alignat}
  where $P_{\mathrm{atm}}$ denotes the constant atmospheric pressure.
\end{subequations}

We restrict our attention to unidirectional waves satisfying
\begin{equation} \label{uni}
  u > 0
\end{equation}
in $D$, and work in dimensionless units where the relative mass flux
\begin{align*}
  m = \int_0^\eta u\, \mathrm dy 
\end{align*}
and gravitational constant $g$ are both normalized to $1$. Lastly, we assume that the vorticity 
\[
\omega = v_x - u_y = -\gamma < 0
\]
is constant and adverse.

Introducing a stream function $\psi$ satisfying $\psi_x = -v$ and $\psi_y = u$, with the normalization that $\psi = 0$ at the bottom, we can reformulate \eqref{eqn:trav} as
\begin{subequations}\label{sys:stream}
  \begin{alignat}{2}
    \label{sys:stream:lap}
    \Delta\psi - \gamma &=0 &\qquad& \text{in } D,\\
    \label{sys:stream:bern}
    \tfrac 12\abs{\nabla\psi}^2 + y  &= r &\quad& \text{on }y=\eta,\\
    \label{sys:stream:kintop} 
    \psi  &= 1 &\quad& \text{on }y=\eta,\\
    \label{sys:stream:kinbot} 
    \psi  &= 0 &\quad& \text{on }y=0,
  \end{alignat}
  where here $r$ is a (non-dimensional) Bernoulli constant. The assumption \eqref{uni}  translates to
  \begin{equation}\label{unipsi}
    \psi_y > 0
  \end{equation}
  everywhere in $D$. 
  
\end{subequations}

\section{Preliminaries} \label{s:preliminaries}

An important part of the analysis is related to laminar flows subject to \eqref{uni}, whose stream functions can be obtained as $\psi = \Psi(y;s)$ where
\[
\Psi(y;s) = \tfrac12 \gamma y^2 + s y
\]
and the parameter $s > 0$ is the horizontal velocity at the bottom. Note that we can restrict to positive values of $s$ because of our assumption \eqref{unipsi}. The corresponding depth and the Bernoulli constant are given explicitly by
\[
  d(s) = \frac{-s + \sqrt{s^2 + 2\gamma}}{\gamma}, 
  \qquad 
  R(s) = \frac12 s^2 + \gamma + d(s).
\]
Thus, for each $s>0$ the pair $\psi = \Psi(y;s)$ and $\eta = d(s)$ satisfies \eqref{uni} and solves $\eqref{sys:stream}$ with Bernoulli constant $r = R(s)$. 
% For more on unidirectional flows we refer to \cite{Kozlov2015}. 

From the above formulas we see that the depth $d(s)$ is a monotonically decreasing function of $s > 0$, while $R(s)$ decreases for $s \in [0,s_c]$ and increases for $s \geq s_c$ for some critical value $s_c>0$ satisfying $R'(s_c) = 0$. More explicitly, $s_c > 0$ is the unique solution to 
\[
s_c - \frac{1}{\gamma} + \frac{s_c}{\gamma \sqrt{s_c^2+2\gamma}} = 0.
\]
Assuming $\gamma$ is large, we have
\[
s_c = \frac{1}{\gamma} + O(\gamma^{-3/2}).
\]
For the corresponding depth, we find
\[
d(s_c) = \frac{-s_c + \sqrt{s_c^2+2\gamma}}{\gamma} = \sqrt{\frac{2}{\gamma}} - \frac{1}{\gamma^2} + O(\gamma^{-5/2}).
\]
Note also that for any $r \in (R(s_c),R(0)]$ there are exactly two roots $s_-(r)<s_+(r)$ to the equation $R(s) = r$. 

The largest possible value for $d = d(s)$ is attained at $s = 0$ and is given by
\[
d(0) = \sqrt{\frac{2}{\gamma}},
\]
while the corresponding Bernoulli constant is the quantity
\[
R(0) = \gamma + \sqrt{\frac{2}{\gamma}}
\]
appearing on the right hand side of \eqref{eqn:rbound}.
It is known from \cite{Kozlov2015} that an arbitrary non-laminar and unidirectional solution must have a Bernoulli constant $r$ that is strictly greater than $R(s_c)$. Furthermore, unidirectional laminar flows have Bernoulli constants strictly less than $R(0)$.

\section{Proof of the main results}

\begin{lemma} \label{lem:crestmin} 
  For any $\lambda \ge \tfrac{\gamma}{2}$, the minimum of $\psi_y - \lambda y$ restricted to the surface is achieved at the crest.
\end{lemma}
\begin{proof}
  This is essentially proved in \cite[Proposition 3.3]{Kozlov2023}, but we briefly outline the argument here for completeness. Consider the auxiliary function 
  \begin{align*}
    f = \tfrac 12 \abs{\nabla\psi}^2 + y + (\lambda-\gamma)\psi.
  \end{align*}
  Straightforward calculations show that
  \begin{align*}
    \Delta f + af_x + bf_y = \frac{2+(2\lambda-\gamma)(2\psi_y + \lambda\abs{\nabla\psi}^2)}{\abs{\nabla\psi}^2} \ge 0,
  \end{align*}
  where the coefficients $a,b$ are given explicitly by
  \begin{align*}
    a = -2\big(f_x - (2\lambda-\gamma)\psi_y )\big)\abs{\nabla\psi}^{-2},
    \qquad 
    b = -2\bigl(-2 + f_y - (2\lambda-\gamma)\psi_y\big)\abs{\nabla\psi}^{-2}.
  \end{align*}
  Moreover, $f$ is constant on the surface, and satisfies $f_y = 1+\lambda \psi_y > 0$ on the bottom. By the maximum principle, $f$ therefore achieves its maximum at every point along the surface. The Hopf lemma then gives 
  \begin{align*}
    0 < f_y = \frac{\eta_{xx}\psi_y^2 + 1}{1+\eta_x^2} + \lambda \psi_y 
    = -\frac{\psi_y}{\eta_x} \frac d{dx} (\psi_y(x,\eta(x)) - \lambda \eta(x))
  \end{align*}
  along the surface. Here the second equality holds wherever $\eta_x \ne 0$, and both equalities are established by repeatedly differentiating the boundary conditions and using the resulting equations to eliminate partials of $\psi$ in favor of derivatives of $\eta$. Thus the total derivatives of $\eta(x)$ and $\psi_y(x,\eta(x))-\lambda\eta(x)$ have opposite signs, and the result follows.
\end{proof}

In what follows the quantity \[\ucrest = \psi_y(0,\hat\eta),\] representing the horizontal relative velocity at the crest, plays a crucial role.

\begin{lemma} \label{lem:etabound} 
  If $\hat \eta \le \frac 2\gamma \ucrest$, then $\hat \eta \le 2/\ucrest$.
\end{lemma}
\begin{proof}
  Consider the harmonic function
  \begin{align*}
    h = \psi_y - \lambda y,
    \qquad 
    \lambda 
    = \frac{\ucrest}{\hat \eta}
    = \frac{\psi_y(0,\hat\eta)}{\hat \eta}.
  \end{align*}
  The choice of $\lambda$ guarantees that $h$ vanishes at the crest. By assumption we also have $\lambda \ge \gamma/2$, so that Lemma~\ref{lem:crestmin} implies that the minimum of $h$ along the surface is achieved at the crest. Thus $h \ge 0$ on the surface. On the bed $h = \psi_y > 0$, and so the maximum principle implies that $h \ge 0$ throughout the fluid. Integrating this inequality along the crest line, we discover
  \begin{align*}
    0 \le \int_0^{\hat \eta} \big(\psi_y(0,y) - \lambda y\big) \, \mathrm dy
    = 1 - \frac \lambda 2 {\hat \eta}^2.
  \end{align*}
  Rearranging, substituting the definition of $\lambda$, and canceling a common factor of $\hat \eta > 0$ yields the desired upper bound.
\end{proof}

\begin{figure}
  \centering
  \includegraphics[scale=1.0]{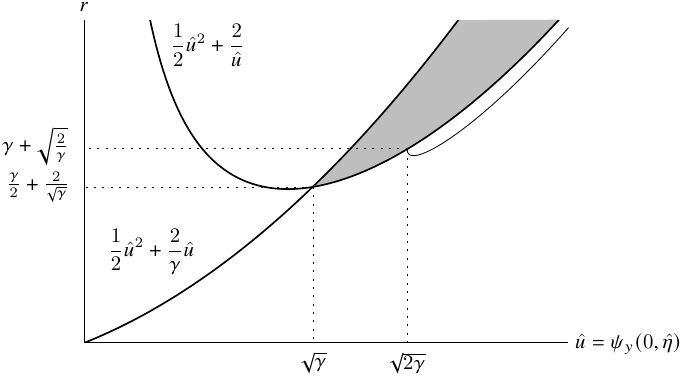}
  \caption{The inequality in \eqref{eqn:nor} is shown in the shaded region bounded by the two thicker curves. The thin curve shows $(\hat u,r)$ values for laminar flows, with a minimum $r$-value of $R(s_c)$. Simple numerical calculations reveal that $R(s_c) > \gamma/2+2/\sqrt{\gamma}$ for $\gamma \ge 1.372$.}
  \label{fig:nor}
\end{figure}

\begin{corollary}\label{cor:r}
  There are no solutions whose Bernoulli constant $r$ lies in the interval
  \begin{align}
    \label{eqn:nor}
    \frac 12 \ucrest^2 + \frac 2\ucrest < r \le \frac 12 \ucrest^2 + \frac 2\gamma \ucrest.
  \end{align}
\end{corollary}
\begin{proof}
  Evaluating the dynamic boundary condition \eqref{sys:stream:bern} at the crest gives $r = \frac 12 \ucrest^2 + \hat\eta$. Thus \eqref{eqn:nor} is equivalent to
  \begin{align*}
    \frac 2\ucrest < \hat \eta \le \frac 2\gamma \ucrest,
  \end{align*}
  or in other words that the hypothesis of Lemma~\ref{lem:etabound} holds but the conclusion does not hold.
\end{proof}

The inequality \eqref{eqn:nor} is sketched as the shaded region in Figure~\ref{fig:nor}.

\begin{proposition}\label{prop:bound}
  Let $\gamma > 0$ be such that $R(s_c) > \gamma/2+2/\sqrt{\gamma}$. Then along any global bifurcation curve of Stokes or solitary wave solutions, connected to a laminar flow, we have
  \begin{align}
    \label{eqn:propbound}
    \psi_y(0,\hat\eta) = \ucrest \ge u_*(r) > \sqrt\gamma,
  \end{align}
  where $u_*(r)$ is the largest root of $u_*^2/2 + 2/u_* = r$. 
\end{proposition}
\begin{proof}
  Any bifurcation curve of Stokes or solitary waves starts with a laminar flow
  located to the right of the shaded region in Figure~\ref{fig:nor}. Recall
  from Section~\ref{s:preliminaries} that any non-laminar solution must have $r >
  R(s_c)$. Thus if $\gamma > 0$ satisfies $R(s_c) > \gamma/2+2/\sqrt{\gamma}$,
  all solutions along the curve will satisfy $r > \gamma/2+2/\sqrt{\gamma}$,
  and hence correspond to points in Figure~\ref{fig:nor} lying above the corner
  of the shaded region. By Corollary~\ref{cor:r} they must also stay to
  the right of the shaded region, and the desired inequality follows.
\end{proof}

For large $r$, the function $u_*(r)$ appearing in \eqref{eqn:propbound} has the expansion
\begin{align}\label{uasymp}
  u_*(r) = \sqrt{2r} - \frac 1r - \frac 3{2^{3/2}} \frac 1{r^{5/2}} + O(r^{-4}).
\end{align}
Moreover, as shown in Figure~\ref{fig:nor}, we have
\begin{align*}
  u_*(R(0)) = \sqrt{2\gamma}.
\end{align*}
Also, as all solutions automatically have
\begin{align}\label{Rcasymp}
  r > R(s_c) = \gamma + \sqrt {\frac 2{\gamma}} + \frac{1}{2\gamma^2}+ O(\gamma^{-5/2}),
\end{align}
the hypothesis $R(s_c) > \gamma/2+2/\sqrt{\gamma}$ holds provided $\gamma$ is sufficiently large. Numerically calculating $R(s_c)$ and plotting it as a function of $\gamma$, it appears that $\gamma \ge 1.372$ is sufficient.

\subsection{Proof of Theorem \ref{thm:main}}

\begin{proof}
  First we establish that there exists exactly one $\gamma = \gamma_\star \approx 1.372$ such that
  \[
  R(s_c(\gamma); \gamma) = \frac{\gamma}{2} + \frac {2}{\sqrt{\gamma}}.
  \]
  For this purpose we consider the function
  \[
  \rho(\gamma) = R(s_c(\gamma), \gamma) - \frac{\gamma}{2} - \frac {2}{\sqrt{\gamma}}.
  \]
  Computing the derivative and taking into account that 
  \[
  \frac{\partial R}{\partial s}(s_c(\gamma), \gamma) = 0,
  \]
  we find
  \[
  \begin{split}
    \rho'(\gamma) & = \frac{\partial R}{\partial \gamma}(s_c, \gamma) - \frac12 + \frac{1}{\sqrt{\gamma^3}} = - \frac{d(s_c)}{\gamma}+\frac12 + \frac{1}{\gamma \sqrt{s_c^2+2\gamma}} + \frac{1}{\sqrt{\gamma^3}}  \\
    & = - \frac{2}{\gamma (s_c + \sqrt{s_c^2+2\gamma})}+\frac12 + \frac{1}{\gamma \sqrt{s_c^2+2\gamma}} + \frac{1}{\sqrt{\gamma^3}} > \frac12 -\frac{1}{\gamma \sqrt{s_c^2+2\gamma}} + \frac{1}{\sqrt{\gamma^3}} > 0.
  \end{split}
  \]
  Thus, since $\rho(\gamma)$ is strictly increasing, negative for small $\gamma$, and is positive for large $\gamma$, there exists a unique zero. A simple numerical computation reveals that $\gamma_\star \approx 1.37184$.
  
  Now consider a bifurcation curve of solutions with $\gamma > \gamma_\star$ starting at a laminar flow. By Proposition~\ref{prop:bound} and the above argument, any solution along this curve satisfies \eqref{eqn:propbound}. From Lemma~\ref{lem:crestmin} we know that $\psi_y - \tfrac12 \gamma y$ has minimum along the surface at the crest, and hence
  \[
  \psi_y - \tfrac12\gamma \eta \geq \psi_y(0,\hat{\eta}) - \tfrac12 \gamma\hat{\eta}
  \]
  everywhere along the surface. From \eqref{eqn:propbound} and the definition of $u_*(r)$ we have $\hat\eta < 2/u_*(r)$, and so the above inequality implies
  \[
    \psi_y \geq \frac12\gamma \eta + u_*(r) - \frac{\gamma}{u_*(r)} > \frac12\gamma \eta,
  \]
  where in the last step we have used $u_*(r) > \sqrt\gamma$. The proof is complete.
\end{proof}

\subsection{Proof of Theorem \ref{thm:main1}}

Let us prove the first claim. Note that by Lemma \ref{lem:crestmin} the function $\psi_y - \gamma y$ attains its minimum along the surface at the crest, which gives
\[
\psi_y \geq \gamma \eta + \hat{u}- \gamma\hat{\eta} \geq \gamma \check{\eta} + \hat{u}- \gamma\hat{\eta}.
\]
Now we need to compute an asymptotic expansion in $\gamma$ for $\hat{u}- \gamma\hat{\eta}$. We rewrite
\begin{align*}
  \hat{u}- \gamma\hat{\eta} 
  &= 
  \hat{u}- \gamma\big(r - \tfrac 12 \hat{u}^2\big) 
  \ge u_*(r) - \gamma\big(r - \tfrac 12 u_*^2(r)\big) \\
  &= 
  \sqrt{2r} - \frac 1r - \frac 3{2^{3/2}} \frac 1{r^{5/2}} + O(r^{-4}) 
  \\&\qquad 
  - \gamma \left( r  - \frac12 \left[\sqrt{2r} - \frac 1r - \frac 3{2^{3/2}} \frac 1{r^{5/2}} 
  + O(r^{-4})\right]^2\right).
\end{align*}
Here we have used Proposition~\ref{prop:bound} and \eqref{uasymp}.
Next we apply \eqref{Rcasymp} and, after an elementary but long computation, obtain $\hat{u}- \gamma\hat{\eta} \geq O(\gamma^{-5/2})$. 
Then it is left to estimate $\check{\eta}$. Using the assumption $r \le R(0)$ we get the upper bound $\check{\eta} < r \le R(0) = \sqrt{2/\gamma}$. For the lower bound we use the inequality $\check{\eta}>d(s_1)$, which follows from the assumption on $r$, and where $s_1 > 0$ is the unique solution to $R(s_1) = R(0)$; see \cite{Kozlov2015}. A simple Taylor expansion shows that $d(s_1) = \sqrt{2/\gamma} + O(\gamma^{-2})$. 

The two remaining statements can be obtained in a similar fashion.

\subsection{Proof of Theorem \ref{thm:main2}}

We shall consider for a fixed $\gamma > \gamma_\star$ the bifurcation curve of solitary waves $\mathcal C$ constructed by \cite{Haziot2023}, containing small-amplitude unidirectional solutions. Let ${\mathcal C}_\uni$ be a connected component of unidirectional waves. First, we note that the Bernoulli constant remains bounded along $\mathcal C_\uni$, which follows from \cite{Lokharu2021}. Thus, it follows from Theorem~\ref{thm:main} that the surface is free from stagnation points. We claim that the norms of solutions are therefore uniformly bounded along $\mathcal C_\uni$. To see this, first note that the fluid depth is uniformly bounded from above and is separated from zero from below along $\mathcal C_\uni$, by using an upper bound for $r$ and bounds from \cite{Kozlov2015}. Then one needs to know that $\psi_y > \delta > 0$ in a neighbourhood of the surface. But this is true by the following reasoning. The function $\psi_y - \varepsilon y$ is harmonic and we can choose $\varepsilon$ (depending only on $\gamma$) so that it is positive at the surface (note that $\psi_y$ is uniformly separated from zero for all unidirectional solutions connected to the laminar flow). Then $\psi_y - \varepsilon y > 0$ everywhere by the maximum principle. 
Now a local norm estimate as in Proposition 3.2 of \cite{Kozlov2023} is valid along a uniform neighbourhood around the surface. Near the bottom we always have a uniform estimate in terms of $\gamma$ by the classical elliptic theory. Now it follows directly by Theorem~1.1 in \cite{Haziot2023} that the bifurcation curve $\mathcal C$ must leave the unidirectional regime. By the maximum principle and a compactness argument, this transition must occur at a smooth solution where $\psi_y$ vanishes somewhere along the bottom. Finally, by applying the Hopf lemma to $\psi_x$ we know that $\psi_y$ attains its minimum along the bottom at the point directly below the crest. The proof is complete.

\noindent{\bf Conflict of interests.} The authors report no conflict of interests. \\

\noindent{\bf Data availability.} We do not analyse or generate any datasets, because our work proceeds within a theoretical and mathematical approach.

\bibliographystyle{sn-bibliography}
\bibliography{bibliography}
\end{document}